\documentclass[11pt]{amsart}
\input epsf.tex
\usepackage{latexsym}
\usepackage{amsmath,amsthm,amsfonts,amssymb}
\usepackage{epsfig}

\theoremstyle{theorem}
\newtheorem{pro}{Proposition}[section]
\newtheorem{thm}[pro]{Theorem}
\newtheorem{lem}[pro]{Lemma}

\newtheorem{cor}[pro]{Corollary}

\theoremstyle{definition}
\newtheorem{defi}[pro]{Definition}

\newtheorem{cnj}[pro]{Conjecture}
\newtheorem{rem}[pro]{Remark}

\def\H{{\mathbb H}}
\def\N{{\mathbb N}}
\def\R{{\mathbb R}}
\newcommand{\titap}{\operatorname{titap}}
\newcommand{\vol}{\operatorname{vol}}
\newcommand{\snap}{\operatorname{\it snap}}
\newcommand{\smear}{\operatorname{\it smear}}

\begin{document}

\title{The Gromov Norm of the Product of Two Surfaces}
\thanks{This paper has been published in \emph{Topology}, {\bf 44}:2 (2005), 321--339.
However, the printed version contains an error that invalidates one direction of the main theorem. The present version contains an erratum, at the end, explaining this}

\author[L.~Bowen]{Lewis Bowen}
\address{Lewis Bowen, Department of Mathematics, 
University of Hawaii, 
 2565 McCarthy Mall, Honolulu, HI 96822, USA}
\email{lpbowen@math.hawaii.edu}

\author[J.~A.~De Loera]{Jes\'us A. De Loera}
\address{Jes\'us A. De Loera, 
Department of Mathematics,
University of California, 
Davis, CA 95616, USA}
\email{deloera@math.ucdavis.edu}

\author[M.~Develin]{Mike Develin}
\address{Mike Develin, 
American Institute of Mathematics, 
360 Portage Ave., 
Palo Alto, CA 94306, USA}
\email{develin@post.harvard.edu}
\author[F.~Santos]{Francisco Santos}
\address{
Francisco Santos,
Dept. de Matem\'aticas, Estad\'{\i}stica y Computaci\'on, 
Universidad de Cantabria,
E-39005 Santander, SPAIN}
\email{santosf@unican.es}



\begin{abstract}
We make an estimation of the value of the Gromov norm of the Cartesian
product of two surfaces. Our method uses a connection between 
these norms and the minimal size of triangulations of the products of two polygons. This allows us
to prove that the Gromov norm of this product is between 32 and 52
when both factors have genus 2. The case of arbitrary genera is easy to deduce form this one.
\end{abstract}

\maketitle

\section{Introduction}

Gromov defined the \emph{simplicial volume} (also now known as the
\emph{Gromov norm}) of a closed orientable manifold $M$ 
[Gromov 1982] as the infimum of the
$l^1$-norms of all singular chains representing the top homology
class of $M$. 
It is an invariant  that quantifies the topological 
complexity of $M$. For example, if $\pi_1(M)$ is amenable then its
Gromov norm, denoted $||M||$, vanishes. But if $M$ admits a metric
of negative curvature, then $||M||$ is positive and finite. In fact,
the Gromov-Thurston theorem [Gromov 1982, Thurston] states that
if $M$ is a hyperbolic manifold then
\begin{equation}
||M||=\operatorname{volume}(M)/\operatorname{volume}(s_n),
\end{equation}
where $s_n$ is an ideal  simplex of maximum volume in $\H^n$
(and $\dim(M)=n)$). This result was
then used to give a topological proof of Mostow rigidity for
hyperbolic manifolds.  

More recently, while studying some data
structure problems, Sleator, Tarjan, and Thurston~\cite{sleatoretal} made explicit
computational connection of the Gromov norm to the size of
minimal triangulations of polytopes and balls. 
They used this relation
to compute the exact combinatorial diameter of the associahedron
(or ``Stasheff polytope''),
one of the most important polytopes in combinatorics. In
this paper we continue the inspection of this interrelation between
topology and polyhedral combinatorics. We investigate the Gromov norm of the product of
two surfaces by relating it to triangulations of the product of two polygons.

\begin{figure}
\begin{center}
                \includegraphics[height=1.5in]{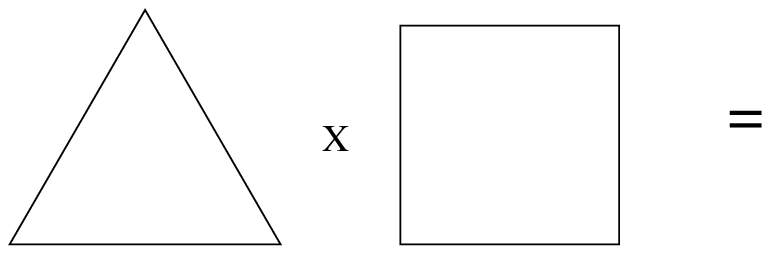}
                \hspace{0.4cm}
                \includegraphics[height=1.5in,width=1.5in]{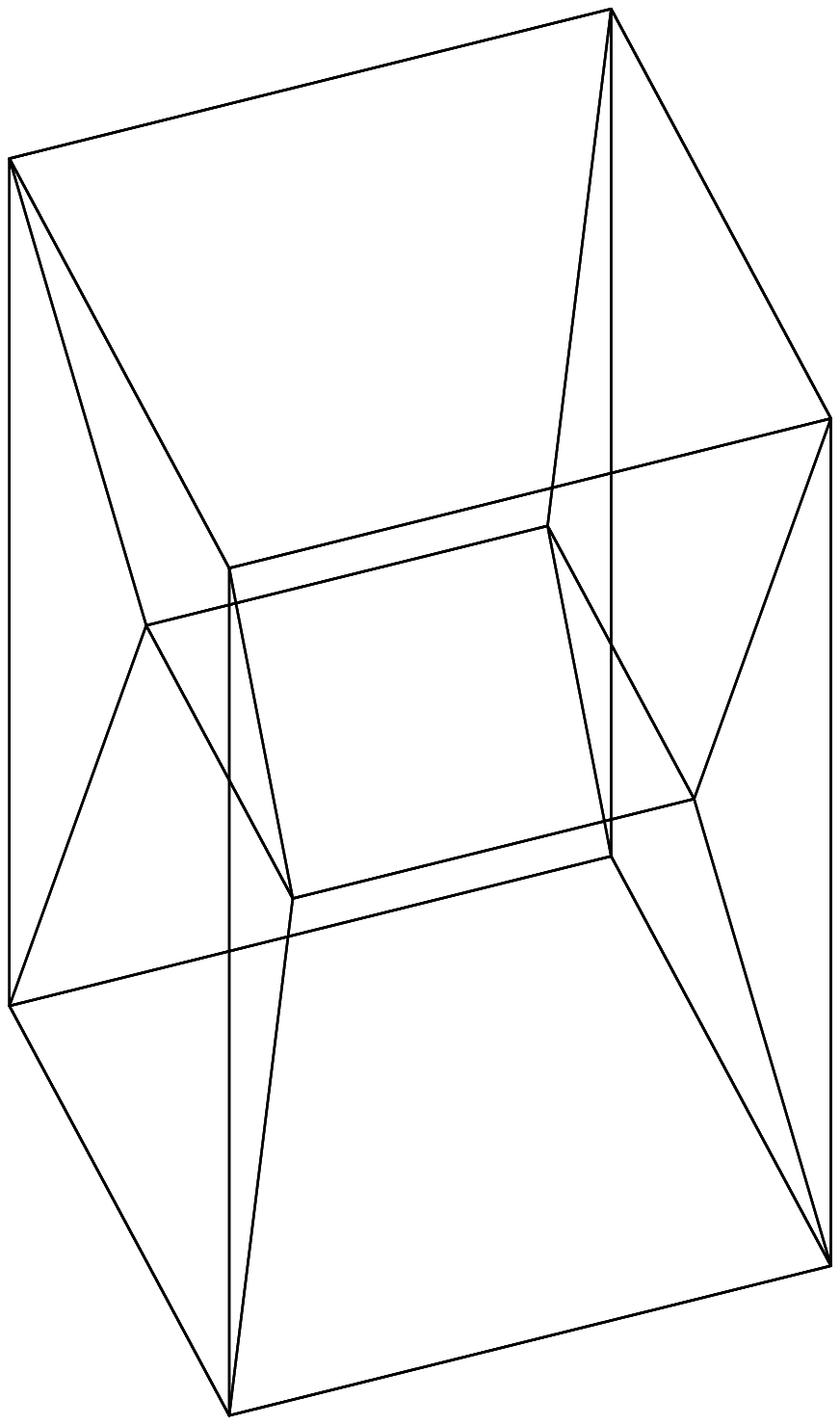}
        \caption{ The Schlegel diagram of the product of a triangle and a square.}
                \label{schlegel}
\end{center}
\end{figure}

More precisely, in Section \ref{polytopalgromov} we define the
 \emph{polytopal Gromov norm} of a convex polytope,
and show that the Gromov norm of the product of two surfaces can be 
computed from the polytopal Gromov norm of the product 
$P(n,m)$ of an $n$-gon with an $m$-gon, for $n$ and $m$ 
asymptotically big. $P(n,m)$ is a four-dimensional polytope with
$m+n$ facets: $m$ prisms over an $n$-gon and $n$ prisms over an
$m$-gon.  In Figure \ref{schlegel} we present the Schlegel diagram of
$P(3,4)$. 

\begin{thm}
\label{surfaces}
Let $||P||$ denote the polytopal Gromov norm of a polytope $P$.
Then, the Gromov norm of the product $\Sigma_g\times\Sigma_h$
of two surfaces of genera $g$ and $h$ equals
\begin{equation}
\frac{||\Sigma_g\times\Sigma_h||}{(g-1)(h-1)}=16\lim_{n,m\to\infty}\frac{||P(n,m)||}{nm}= 16\inf_{n,m}\frac{||P(n,m)||}{nm}. 
\label{eqn-covering}
\end{equation}
\end{thm}

The case $g=h=2$ of Theorem \ref{surfaces} is proved in Section~\ref{polytopalgromov}. The 
general case follows from the following well-known lemma, since  $\Sigma_g\times\Sigma_h$ is a $(g-1)(h-1)$-fold covering
of $\Sigma_2\times\Sigma_2$.

\begin{lem} \label{Gromov}
If $f:M \to N$ is a degree $\deg(f)$ map between closed orientable
manifolds $M$ and $N$ then
\begin{equation}
||M|| \ge \deg(f)||N||
\end{equation}
If $f$ is a covering map then the above inequality is an equality.
\end{lem}

The
polytopal Gromov norm of a polytope $P$ is, roughly speaking, the minimal cardinality of
an affine triangulation of $P$ ``with real coefficients''
 (see the precise definition in
Section~\ref{polytopalgromov}). By definition, it is at most
equal to the minimum number of simplices needed to
(affinely) triangulate $P$. For this reason, and for its intrinsic
interest, we try in section \ref{triangus} to compute the size of
a minimal triangulation of $P(m,n)$. Our main results in this 
direction are:

\begin{thm}   \label{triangulations}
Let $T(m,n)$ denote the minimal number of simplices in a triangulation of
$P(n,m)$. Then:
\begin{equation}
\label{eq:3xm}
 \forall m \text { odd: }\quad
\left\lfloor \frac{9m-15}{2} \right\rfloor \  \le \  ||P(3,m)||\  \le  \  \left\lceil \frac{9m-15}{2} \right\rceil \ =\ T(3,m),
\end{equation}
\begin{equation}
 \forall m \text { even: }\quad
\frac{15}{2}(m-2)\   \le  \  ||P(4,m)||\  \le \ T(4,m)\ \le \  8(m-2), \\
\end{equation}
\begin{equation}
\forall m, n \text { even: }\qquad
 T(m,n)\ \le \  \frac{7}{2}mn - 6(m+n) +8.
\label{sevenhalves}
\end{equation}

\end{thm}

In Section 3 some of these statements are more detailed and do not require any parity condition on $m$ and $n$.
For example, for even $m$ we know that $9m/2-9 \le ||P(m,3)|| \le 9m/2-8 =T(3,m)$. Observe, however, 
that in order to apply Theorem \ref{surfaces} we only need to know $||P(m,n)||$ for a sequence of values
of $(m,n)$ with both $m$ and $n$ going to infinity. In particular, 
the result in equation~(\ref{sevenhalves}) already gives an upper bound for the Gromov norms we are interested in.
In Section~\ref{upperbound} we get a better upper bound by constructing a {\em binary cover} with asymptotically fewer simplices,
in the case $m=n$. A general lower bound for $P(m,n)$ is computed in Section~\ref{lowerbound}:

\begin{thm} \label{mainthm}
\begin{eqnarray}
||P(m,m)|| &\le& \frac{13}{4}m^2-\frac{19}{2}m, \quad \text{ for even $m$},\\
||P(n,m)|| &\ge& 2mn-\frac{8}{3}(m+n).
\end{eqnarray}
\end{thm}

Theorems \ref{mainthm} and \ref{surfaces} together give the following corollary, which is our main
result: 

\begin{cor}
For every positive integers $g$ and $h$, the Gromov norm of the product of $\Sigma_g$ and
$\Sigma_h$ is bounded by:
    $ 32 \leq \frac{||\Sigma_g \times \Sigma_h||}{(g-1)(h-1)} \leq 52$.
\end{cor}

Let us compare this result with earlier bounds. It is well-known that
the Gromov norm of a genus $g$ surface is zero if $g \le 1$, and is
equal to $4(g-1)$ otherwise. Hoster and Kotschick \cite{hok} proved
that whenever $M$ is a surface bundle over a surface $B$,
with fiber $||F||$, then $||M|| \ge
||F||||B||$. This implies that $\frac{||\Sigma_g \times
\Sigma_h||}{(g-1)(h-1)} \ge 16$. In contrast our new lower
bound is 32. On the other end, from any triangulation of the surface $\Sigma_2$ with $T$
triangles, the product tiling of $\Sigma_2 \times \Sigma_2$ via the
Cartesian product of two triangles can be subdivided into $6T^2$
triangles. This gives the easy (and known) result $\frac{||\Sigma_g
\times \Sigma_h||}{(g-1)(h-1)} \le 96$, but now we have 52 as an
upper bound.


\section{Polytopal Gromov norm}\label{polytopalgromov}

We start by recalling the detailed definition of the Gromov norm of
a closed orientable manifold $||M||$. Let $S(M)$ be the singular chain complex of $M$ (with
real coefficients). For each nonnegative integer $k$, $S_k(M)$ is a
real vector space with basis consisting of all continuous maps
$\sigma:\Delta^k \to M$ where $\Delta^k$ is a
$k$-simplex. The norm of an element

\begin{equation} \label{homologelem}
c=\Sigma_\sigma \, r_\sigma, \qquad \sigma \in S_k(M) 
\end{equation} 
is defined by
\begin{equation}
||c||:=\Sigma_\sigma |r_\sigma|.
\end{equation}

If $\alpha$ is a homology class in $H_k(M)$, its simplicial norm is by
definition the infimum of $||c||$ over all $k-$chains $c \in S_k(M)$
representing $\alpha$. The
\emph{Gromov norm} of $M$ is the simplicial norm of the fundamental
class $[M] \in H_n(M)$.

In this section we relate the Gromov norm of $\Sigma_2\times \Sigma_2$ with
what we call  the \emph{polytopal Gromov norm} of $P(n,m)$:

\begin{defi} Let $P$ be a polytope. For each $k\in\N$,
let $S_k(P)$ be the $\R$-vector space with basis equal to the set of
all \emph{affine} maps $\sigma: \Delta^k \to P$.  We call a chain in
$S_k(P)$ an affine chain. We let $S(M)$ denote the resulting singular
chain complex and
\begin{equation}
S(P,\partial P)=S(M)/S(\partial M)
\end{equation}
be the relative chain complex of $(P, \partial P)$. We denote the
resulting homology by $H_*(P, \partial P)$. As before, we define a
norm on the chains in $S(P)$ which induces a pseudonorm on the chains
of $S(P, \partial P)$. In turn, this induces a pseudonorm on $H_*(P, \partial
P)$. The \emph{polytopal Gromov norm} (or \emph{polytopal simplicial
volume}) of $P$ is the pseudonorm of the fundamental class of $[P,
\partial P] \in H_n(P, \partial P)$. We denote it by $||P||$.
\end{defi}

\begin{rem}
\label{simplechains}
 If we do not require $\sigma$ to be affine, then we
would be left with the usual definition of the relative Gromov
norm. But the relative Gromov norm of $P$ (with respect to its
boundary) is zero since it is homeomorphic to a ball, which admits
self-maps of arbitrary degree. Instead of requiring that each map
$\sigma$ is affine, we could require only that $\sigma$ takes every
face of $\Delta^k$ into a single face of $P$. The resulting norm gives
the same value since any such map can be ``straightened'' into a
affine map that agrees on the vertices. 
\end{rem}

\begin{thm}
\begin{equation}
||\Sigma_2 \times \Sigma_2|| = \lim_{n,m \to \infty} \,  \frac{16||P(n,m)||}{(n-2)(m-2)}.
\end{equation}
\end{thm}

We prove this in the next four lemmas.

\begin{lem}\label{inflem}
\begin{equation} \label{limit1}
\lim_{n,m \to \infty} \, \frac{||P(n,m)||}{(n-2)(m-2)} = \inf_{n,m} \frac{||P(n,m)||}{(n-2)(m-2)}.
\end{equation}
\end{lem}

\begin{proof}
It suffices to show that for every fixed positive integers  $j$ and $k$
and for sufficiently large $m$ and $n$,

\begin{equation}
\frac{||P(j,k)||}{(j-2)(k-2)} \ge \frac{||P(n,m)||}{(n-2)(m-2)} +{\rm O}((j+k)/mn).
\label{limit2}
\end{equation}
 
To start, suppose that $j-2$ divides $n-2$
and $k-2$ divides $m-2$. Divide the $n$-gon into $\frac{n-2}{j-2}$
$j$-gons and the $m$-gon into $\frac{m-2}{k-2}$ $k$-gons. Taking the product, we obtain a
partition of $P(n,m)$ into $\frac{(n-2)(m-2)}{(j-2)(k-2)}$ copies of
$P(j,k)$.
\begin{figure}
  \centerline{
   \includegraphics[height=1.75in]{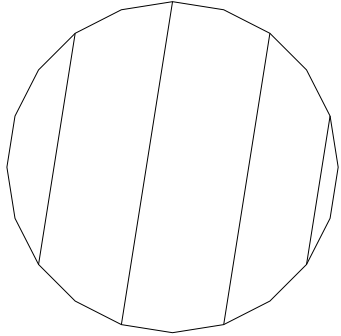}                }
   \caption{ $20$-gon chopped into four hexagons and a quadrangle.} \label{choping}
\end{figure}
From any chain $c$ of $S(P(j,k), \partial P(j,k))$ representing the
fundamental class we can construct a fundamental chain ${\tilde c}$ on
$S(P(n,m), \partial P(n,m))$ by
combinatorially reflecting the chain $c$
of any particular copy of $P(j,k)$ to the adjacent ones. The new chain
satisfies
\begin{equation}
||{\tilde c}|| = ||c||\frac{(n-2)(m-2)}{(j-2)(k-2)}.
\end{equation}
Equation (\ref{limit2}) above 
follows in this case by choosing a sequence of chains $c_i$ so that $||c_i|| \to
||P(j,k)||$.

If $n$ and $m$ are large but $j-2$ does not divide $n-2$
and/or $k-2$ does not divide $m-2$, let $n'$ and $m'$
be the first integers after $n$ and $m$ and such that 
$j-2$ divides $n'-2$ and $k-2$ divides $m'-2$. Clearly
\[
n\le n' < n+j-2,\qquad 
m\le m' < m+k-2.
\]
%
Since $||P(a,b)||$ is an increasing function of $a$ and $b$ (because we can always collapse
chains from bigger to smaller polytopes), we conclude that
\begin{equation}
||P(n,m)|| \le ||P(n',m')|| \le  \frac{(n+j-4)(m+k-4)}{(j-2)(k-2)}||P(j,k)||,
\end{equation}
from which equation (\ref{limit2}) follows.
\end{proof}

In the next step we need the following standard results from 
hyperbolic geometry. See \cite{Rat}.

\begin{lem}\label{hyperbolic}
For every pair of positive integers $n,j$ such that $(n-2) - n/j > 0$
there is a regular geodesic $n$-gon in the hyperbolic plane with
interior angles equal to $\pi/j$. The area of this $n$-gon is $((n-2)
- n/j)\pi$. The group $G_{n,j}$ generated by reflections in its sides
is discrete and acts cocompactly on the hyperbolic plane. There is a
torsion-free finite index subgroup $G'_{n,j} < G_{n,j}$ consisting of
orientation-preserving isometries. 
\end{lem}

\begin{lem}
For any pair of integers $n,m \ge 3$, 
\begin{equation}
||\Sigma_2 \times \Sigma_2|| \le \frac{16}{(n-2)(m-2)}||P(n,m)||.
\end{equation}
\end{lem}
\begin{proof}

Let $c$ be any chain in $S(P(n,m), \partial P(n,m))$ and suppose that $c$ represents $[P(n,m), \partial P(n,m)]$. Choose $j,k$ large enough so that

\begin{eqnarray}
&& (n-2) - n/j >0 \hbox{ and} \\
&& (m-2) - m/k >0.
\end{eqnarray}

We represent $P(n,m)$ as the polytope in $\H^2 \times \H^2$ formed
from a regular $n$-gon with interior angles equal to $\pi/j$ times a
regular $m$-gon with interior angles equal to $\pi/k$. Because the
$n$-gon and the $m$-gon both tile $\H^2$ by reflection, $P(n,m)$ tiles 
$\H^2 \times \H^2$ by
reflections. Using these reflections, the chain $c$
induces a chain ${\tilde c}$ on $\H^2 \times \H^2$
($\H^2$  denotes the hyperbolic plane)  that is invariant
under $G:=G_{n,j} \times G_{m,k}$. But this group has a finite index
subgroup $G':=G'_{n,j} \times G'_{m,k}$ that is torsion-free with no
orientation-reversing elements. Hence the chain ${\tilde c}$ pushes
forward to a chain $c_M$ on the quotient manifold $M$. 
Let $c_M$ represent the fundamental
class $[M]$.  From this construction, we have

\begin{equation}
||c_M|| = [G:G']||c||.
\end{equation}

It is easy to see that $M$ is equal to the Cartesian product of two surfaces with both
surfaces orientable and of genus at least two. Therefore, there is a covering
map $\pi: M \to \Sigma_2 \times \Sigma_2$. Let $c_\Sigma =
\pi_*([c_M])/deg(\pi)$. It represents the fundamental class $[\Sigma_2
\times \Sigma_2]$. Its norm satisfies $||c_\Sigma|| =
||c_M||/deg(\pi)$. So,

\begin{eqnarray}
||c_\Sigma|| &=& ||c_M||/deg(\pi)\\
        &=& ||c||[G:G']/deg(\pi)\\
        &=& ||c|| \vol(\Sigma_2 \times \Sigma_2)/\vol(P(n,m))\\
        &=& ||c|| \frac{16\pi^2}{(n-2 -n/j)(m-2-m/k)\pi^2}.
\end{eqnarray}
The third equality above follows by observing that the volume of $M$ is equal to $[G:G']\,\vol(P(n,m))$ since $M$ is
tiled by $[G:G']$ copies of $P(n,m)$. Similarly, the volume of $M$ is equal to $deg(\pi)$ times the volume of
$\Sigma_2 \times \Sigma_2$. Letting $j$ and $k$ tend to infinity in the above and letting $||c||$ tend towards
$||P(n,m)||$ finishes the lemma.

\end{proof}

\begin{lem}
\begin{displaymath}
||\Sigma_2 \times \Sigma_2|| \ge \inf_{n,m} \, \frac{16}{(n-2)(m-2)}||P(n,m)||.
\end{displaymath}
\end{lem}
\begin{proof}

Let $G$ be a discrete group of hyperbolic isometries such that $\H^2/G
\equiv \Sigma_2$. Let $c$ be a chain in $S(\Sigma_2 \times \Sigma_2)$
representing the fundamental class. There is a finite chain $c_0$ of
simplices in $\H^2 \times \H^2$ such that $c = \pi_*(c_0)$ where
$\pi:\H^2 \times \H^2 \to \Sigma \times \Sigma$ is the universal
covering map. It will be convenient to use Thurston's smearing
construction \cite{Thu}. Let
\begin{displaymath}
c_0 = \Sigma_{i=1}^k \, r_i \sigma_i.
\end{displaymath}
For each $i$, we let $\smear(\sigma_i)$ denote the real-valued measure supported on all orientation-preserving isometric translates of $\sigma_i$ that is induced from Haar measure on $Isom^+(\H^2 \times \H^2)$. We define $\smear(c_0)$ by
\begin{displaymath}
\smear(c_0) = \Sigma_{i=1}^k \, \smear(\sigma_i).
\end{displaymath}
If $G'$ is any discrete torsion-free
cocompact group of isometries, $\smear(c_0)$ induces a measure
$\smear_{G'}(c_0)$ on singular simplices of the quotient $(\H^2 \times
\H^2)/G'$ by setting $\smear_{G'}(c_0)(S)=\smear(c_0)({\tilde S})$ where
${\tilde S}$ is any Borel set of singular simplices in $\H^2 \times
\H^2$ that projects down to $S$ injectively. From the construction,
the total mass of $\smear_{G'}(c_0)$ is given by

\begin{equation}\label{volume}
||\smear_{G'}(c_0)|| = ||c||\frac{\vol\left((\H^2 \times \H^2)/G'\right)}{\vol(\Sigma_2 \times \Sigma_2)}.
\end{equation}

The benefit of this formula is that we can pass from the original
chain $c$ defined on $\Sigma_2 \times \Sigma_2$ to a measure-chain on
a more convenient manifold. To define that manifold, let $d$ be the
maximum diameter of the image of $\sigma_i$ in $\H^2 \times \H^2$
($i=1,\dots,k$).  For every $h>0$, there is a regular $4h$-gon $F_h$ with
all interior angles equal to $2\pi/4h$. We choose $h$ large enough so
that every pair of nonadjacent sides of $F_h$ is at least a distance
$d$ apart. We let $T$ be a tiling of the plane $\H^2$ with copies of
$F_h$.

We will need the concept of a \emph{straight simplex}.  We use the Lorenz model of the hyperbolic plane \cite{Rat}. The Lorenz inner product on $\R^3$ is defined by
\begin{displaymath}
x \circ y = -x_1y_1 + x_2y_2 + x_3y_3.
\end{displaymath}
$\H^2$ is identified with the set of vectors $x$ satisfying $x \circ x
= -1$ and $x_1 >0$. We let $|||x|||$ denote the absolute value of
$\sqrt{x \circ x}$. We say that a simplex $\sigma: \Delta^4 \to \H^2$
is straight if for every $x=\Sigma_i \, x^i e_i \in \Delta^4$
\begin{displaymath}
\sigma(x) = \Sigma_i \, x^i\sigma(e_i) / |||\Sigma_i \, x^i \sigma(e_i)|||.
\end{displaymath}
We say that a simplex $\sigma: \Delta^4 \to \H^2 \times \H^2$ is
straight if composing $\sigma$ with a projection to either of the
$\H^2$ factors results in a straight simplex. Such a simplex is
uniquely determined by its vertices. Its image is equal to the convex
hull of its vertices.

Our last operation is called ``snapping''. If $\sigma$ is a singular
simplex in $\H^2 \times \H^2$, we let $\snap(\sigma)$ be the straight
simplex satisfying the following. For all $i=0, \dots, 4$,
$\snap(\sigma)(e_i)$ is the closest vertex of the tiling $T \times T$
to $\sigma(e_i)$ if there is only one closest
vertex and $\snap(\sigma)(e_i)=\sigma(e_i)$ otherwise. The map $\snap$ on
singular simplices induces a map $\snap_*$ on measures on the set of
singular simplices. We will use $\snap_*(\smear(c_0))$ to construct a
chain on $(F \times F, \partial(F \times F))$. But first, if $\sigma$
has diameter at most $d$ then we need to show that $\snap(\sigma)$ is
contained in single tile of $T \times T$.

 To see this, consider the dual tiling $T^*$ of $T$. The vertices of
 $T^*$ are the centers of the tiles of $T$ and there is an edge
 between two dual vertices if and only if the corresponding tiles in $T$ are
 adjacent. In our case, the tiles of $T^*$ are copies of $F$. Since
 every two nonadjacent sides of $F$ have a distance at least $d$
 apart, the image of $\pi_i \sigma$ does not overlap any nonadjacent
 edges of the dual tiling (where $\pi_i$ is projection from $\H^2
 \times \H^2$ onto the $i$-th $\H^2$ factor). Hence there is a vertex
 of the dual tiling contained in all the dual tiles that contain the
 image. Since vertices of the dual correspond to faces of the domain
 tiling $T$ this implies that there is a single tile $\tau$ of $T$
 such that: for every point $x$ in the image of $\pi_i \sigma$, the closest 
 vertex $v$ of $T$ to $x$ is contained in $\tau$. By construction, this
 implies that the projection of $\snap(\sigma)$ to this $\H^2$-factor
 is contained in $\tau$. Since this is true for both $\H^2$ factors,
 $\snap(\sigma)$ is contained in a single tile of $T \times T$.

Let $c_F$ denote the restriction of $\snap_*(\smear(c_0))$ to the set of
simplices that map into a chosen fixed tile of $T \times T$. By
construction, $\snap_*(\smear(c_0))$ is invariant under the symmetries
of tiling $T \times T$ so it is irrelevant which tile we use. $c_F$ is
supported on a finite set of simplices by construction, so we may
identify it with a finite chain. It is a cycle representing $(F \times
F, \partial(F \times F))$ because snapping and smearing commute with
the boundary map. The norm of $c_F$ is
\begin{equation}\label{bigeqn}
||c_F|| = ||\smear_{G'}(c_0)||=\frac{\vol(F \times F)}{\vol(\Sigma_2 \times \Sigma_2)} ||c|| = (h-1)^2 ||c||,
\end{equation}
where $G'$ is a group having $F \times F$ as its fundamental
domain. The first equation holds because the snapping operation
preserves the total mass of a measure under a quotient by a group that
stabilizes the tiling $T \times T$. The second equality is equation
\ref{volume}. The third equality is from Lemma \ref{hyperbolic}.
Since $F \times F$ is a realization of the polytope $P(4h,4h)$,
$||P(4h,4h)|| \le ||c_F||$. Equation \ref{bigeqn} now implies (by
taking $||c||$ arbitrarily close to $||\Sigma_2 \times \Sigma_2||$)
that

\begin{displaymath}
||\Sigma_2 \times \Sigma_2|| \ge \frac{||P(4h,4h)||}{(h-1)^2}.
\end{displaymath}

By Lemma~\ref{inflem}, as $h$ tends to infinity, the right-hand side approaches

\begin{displaymath}
\inf_{n,m} \frac{16 ||P(n,m)||}{(n-2)(m-2)}.
\end{displaymath}

\end{proof}

\section{Small Triangulations of $P(m,n)$} \label{triangus}

Since the size of any triangulation of $P(m,n)$ is an upper bound for its polytopal Gromov norm,
in this section we give bounds, and in some cases exact values, for the size $T(n,m)$
of a minimal triangulation of $P(n,m)$. We consider this optimization
over \emph{all possible coordinatizations} of $P(n,m)$. In this
section we look at the combinatorics of the polytopes $P(n,m)$ and
their triangulations.  We begin with a table of known sizes of minimal
triangulations in specific instances.

\begin{table}[h] 
        \[
        \begin{array}{r|rrrrrrr}
                  & 3 & 4 & 5 & 6 & 7 & 8 & 9 \\\hline
                3 & 6 & 10 & 15 & 19 & 24 & 28 & 33 \\
                4 & 10 & 16 & 26 & 32 & 42 & \le48 & \le58 \\
                5 & 15 & 26 & 38 & \le49 & \le61 & \le72 \\
                6 & 19 & 32 & \le49 & \le60 & \le77 & \le90 \\
                7 & 24 & 42 & \le61 & \le77 \\
                8 & 28 & \le48 & \le72 & \le90 \\
                9 & 33 & \le58 \\
        \end{array}
        \]
        \caption{Minimal size triangulations for $n$-gons times $m$-gons.}
        \label{optitriang}
        \end{table}

For computing Table \ref{optitriang} we followed the approach 
of \cite{dhss}, based on the solution of an integer programming
problem. We think of the triangulations of a polytope as the vertices
of the following high-dimensional polytope: Let $A$ be a
$d$-dimensional polytope with $n$ vertices.  Let $N$ be the number of
$d$-simplices in $A$.  We define $P_A$ as the convex hull of the set
of incidence vectors of all triangulations of ${A}$. For a
triangulation $T$ the {\em incidence vector} $v_T$ has coordinates
$\,(v_T)_\sigma = 1$ if $\sigma \in T$ and $\,(v_T)_\sigma = 0$ if
$\sigma \not\in T$.  The polytope $P_A$ is the {\em universal
polytope} defined in general by Billera, Filliman and Sturmfels
\cite{BIFIST} although it appeared in the case of polygons in
\cite{dantzig}. In \cite{dhss}, it was shown that the vertices of
$P_A$ are precisely the integral points inside a polyhedron that has a
simple description. The rational vertices of this polytope are in
correspondence with the fractional face-to-face covers. The concrete integer
programming problems were solved using {\em C-plex Linear
Solver}$^{TM}$.  The program to generate the linear constraints is a
small $C^{++}$ program written by De Loera and Peterson.

In the rest of the paper we will often use the following result, 
first proved (for triangulations)
in \cite{prismpaper}. The same result (rounded up),
and with almost the same proof, holds for odd $m$ but we do not need it.

\begin{thm}
\label{prism}
Let $m\ge 4$ be an even number.
The minimum triangulation of the prism over an $m$-gon, in any coordinatization,
has size $\frac{5}{2}(m-2)$. This number equals also the polytopal Gromov norm of the prism.
\end{thm}

\begin{proof}
To see that $\frac{5}{2}(m-2)$ is an upper bound for both numbers it suffices to describe
a triangulation of that size. This goes as follows: first, chop alternate vertices of the $m$-prism,
using $m$ tetrahedra, to obtain an $(m/2)$-antiprism. This has $m$ triangular faces and
two polygons of size $m/2$. Triangulate one of them arbitrarily, and triangulate the antiprism by coning
from a vertex in the opposite face. One needs $m-3+m/2-2$ tetrahedra for this.

For the lower bound, we show that every affine chain in the fundamental class has norm at least 
$\frac{5}{2}(m-2)$. Without loss of generality, we assume that the chain has 
all its vertices on vertices of the prism (see Remark \ref{simplechains}). Also, since we
are dealing with homology relative to the boundary, we assume that the chain has no tetrahedron
contained in the boundary. In particular, each tetrahedron is of one of three types:
 a ``bottom tetrahedron''
with three vertices in the bottom  $m$-gon and a vertex in the top, a ``top tetrahedron'' (the converse), or
a ``middle tetrahedron'' with two vertices on each. 
It is obvious that we need at least $m-2$ bottom and $m-2$ top tetrahedra to cover 
the bottom and top $m$-gons. The result then follows if we prove that the number of middle tetrahedra
is at least half that number. This holds because each middle tetrahedron has two
``bottom triangles'' (the ones with two vertices in the bottom) and the projections of these must
also cover the bottom $m$-gon: any vertical (but otherwise generic) line must be covered by a sequence
of tetrahedra in the chain that starts with a bottom tetrahedron and finishes with a top tetrahedron, 
which implies that in between necessarily some middle tetrahedra is used.
\end{proof}

In the case of a triangle times an $m$-gon, the patterns shown in Table
\ref{optitriang} suggested the following result, which is a rephrasing of equation (\ref{eq:3xm})
in Theorem \ref{triangulations}.

\begin{thm}
\label{thm:c3xcm}
In any coordinatization, the minimum-size triangulations 
and the polytopal Gromov norm of $P(3,m)$ satisfy:
\begin{enumerate}
\item If $m$ is odd, $||P(3,m)||=T(3,m)=  9m/2 -15/2 $.
\item If $m$ is even, $T(3,m)= 9m/2-8$ and $||P(3,m)||$ lies
between that number and $9m/2-9$.
\end{enumerate}
\end{thm}

\begin{proof}
Let $C_3$ and $C_m$ denote
the triangle and $m$-gon of which $P(3,m)$ is the product.
Let $A$, $B$ and $C$ denote the three
vertices of $C_3$.  

We first prove the lower bound 
for the norm of an affine simplicial chain (and, hence, for the size of a triangulation).
We assume without loss of generality that all the
vertices in the chain are vertices of $P(3,m)$ and that no 4-simplex
is contained in the boundary of $P(3,m)$.
Then every maximal simplex in the chain falls into one
of the following types:
\begin{enumerate}
\item A ``type $A$'' simplex, with three vertices on $A\times C_m$ and
one in each of $B\times C_m$ and $C\times C_m$. There are at least 
$m-2$ of them (counted with their coefficients) in every affine simplicial chain. 
Similarly, there will be  $m-2$ simplices of types $B$ and $C$.

\item A ``type AB'' simplex, with two vertices on $A\times C_n$, two
on $B\times C_n$ and one on $C\times C_n$. Together with the type
A and type C simplices, these must cover the prism $AB\times
C_m$. Hence, as in the proof of Theorem \ref{prism},
there are at least $\lceil(m-2)/2\rceil$ of them (there can
certainly be more). Similarly, there are at least $\lceil(m-2)/2\rceil$ simplices of
types $AC$ and $BC$.
\end{enumerate}

Adding up these numbers gives $3m +3\lceil m/2\rceil -9$, which coincides
with the stated lower bound for $||P(3,m)||$ in both the odd and even cases.
In the even case, however,  no triangulation with exactly $9m/2-9$
 simplices exists, hence increasing the lower bound by one. This is so because
such a triangulation
 must triangulate each of the three facets $AB\times C_m$,
$AC\times C_m$, and $BC\times C_m$ in its minimal way. But, by the
analysis in \cite{prismpaper}, every minimum-size triangulation of 
an $m$-prism with even $m$ must be obtained 
(as in the proof of Theorem~\ref{prism})
by first cutting alternate
corners and then triangulating the remaining anti-prism. In
particular, the three $m$-gons $A\times C_m$, $B\times C_m$ and
$C\times C_m$ are triangulated by first cutting half the corners, and
the corners cut should be opposite in the three $m$-gons, which is
impossible. This proves $T(3,m)\le 9m/2-8$ in this case.

The proof of the upper bound for $T(3,m)$ is via the explicit construction
of a triangulation with the stated size. The triangulation is depicted 
for $P(3,m)$ in
Figure~\ref{fig:12gon-ABC} in a ``Cayley Trick view''.  The \emph{Cayley
Trick} is a simple but clever construction that, in our case, gives a natural bijection
between the triangulations of $P(3,m)$ and the ``mixed subdivisions'' of
the Minkowski sum of three equal copies of $C_m$ (see \cite{huberetal,San} for
details).

        \begin{figure}[htb]
             \centerline{
                \includegraphics[height=2.7in]{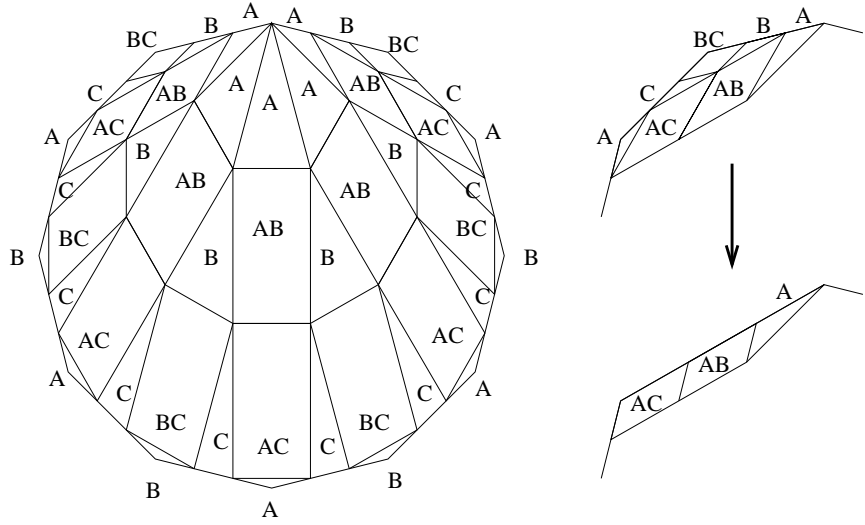}
              }
        \caption{ The minimal triangulation of $C_{12}\times C_3$ (left) and how
to get the one of $C_{11}\times C_3$ from it.}
                \label{fig:12gon-ABC}
        \end{figure}

The triangulation displayed has the number of simplices of types
``A'', ``B'', ``C'', ``AB'' and ``AC'' predicted in the above
paragraphs, and only one more than predicted simplex of type ``AC''.
Exactly the same construction can be done for every even $m$, and produces 
$9m/2-8$ simplices. For odd $m$, we show on the right part of 
the figure how to obtain the minimal triangulation of $P(3,m)$ from that of
$P(3,m+1)$.
%
\end{proof}



The above result suggests a simple and relatively efficient way of
triangulating $C_n\times C_m$: triangulate $C_n$ into $n-2$ triangles
and triangulate each of the resulting $C_3\times C_m$'s in the optimal
way. It is easy to make the triangulations match in common boundaries:
just label the vertices of $C_n$ with $A$, $B$ and $C$ in such a way
that every triangle gets the three labels 
(as in Figure~\ref{fig:Cn})
and replicate the
triangulation of Figure~\ref{fig:12gon-ABC} so that the labels match.
This procedure produces approximately $9mn/2$ maximal simplices.
But this number can be decreased, as follows:

       \begin{figure}[htb]
                \includegraphics[height=1.5in]{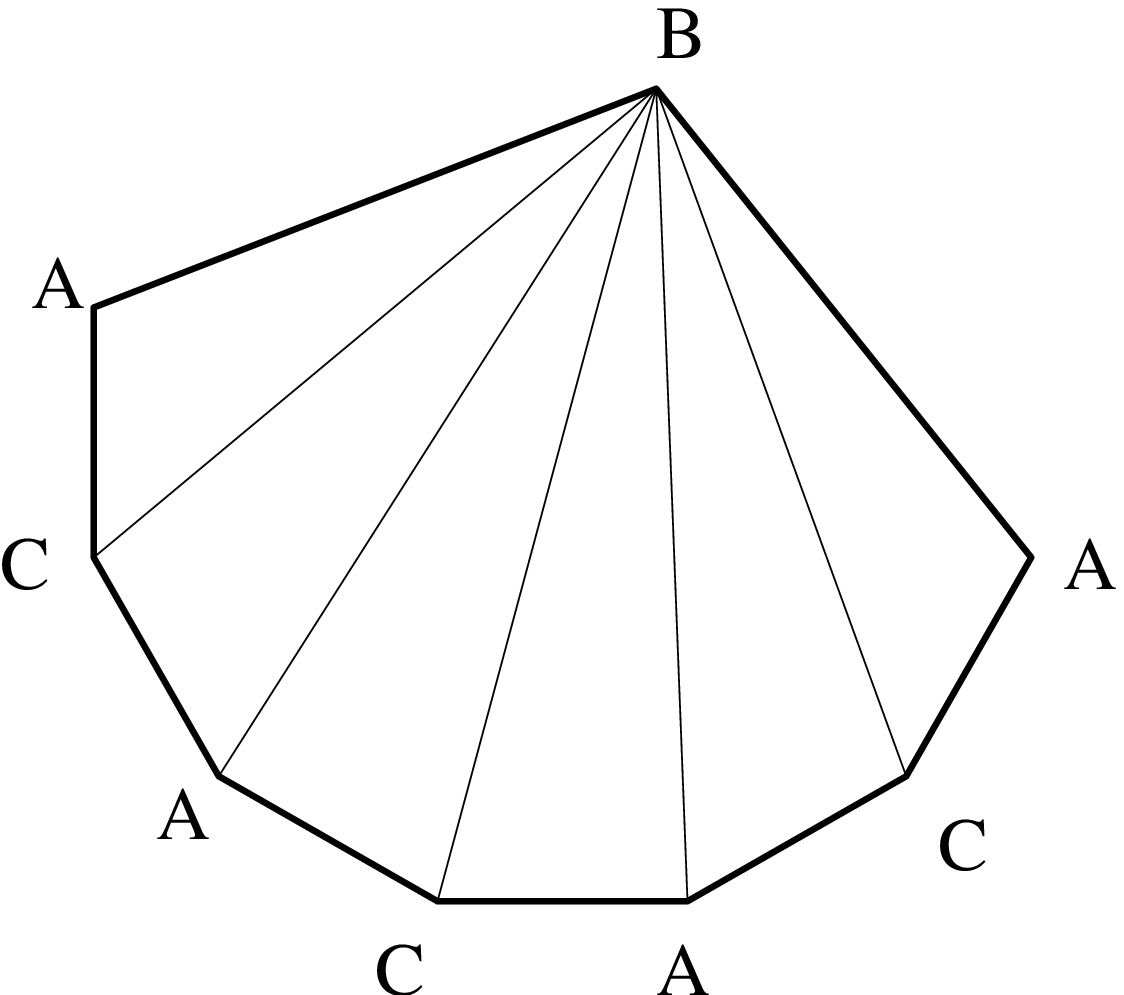}
\hspace{0.5in}
                \includegraphics[height=1in]{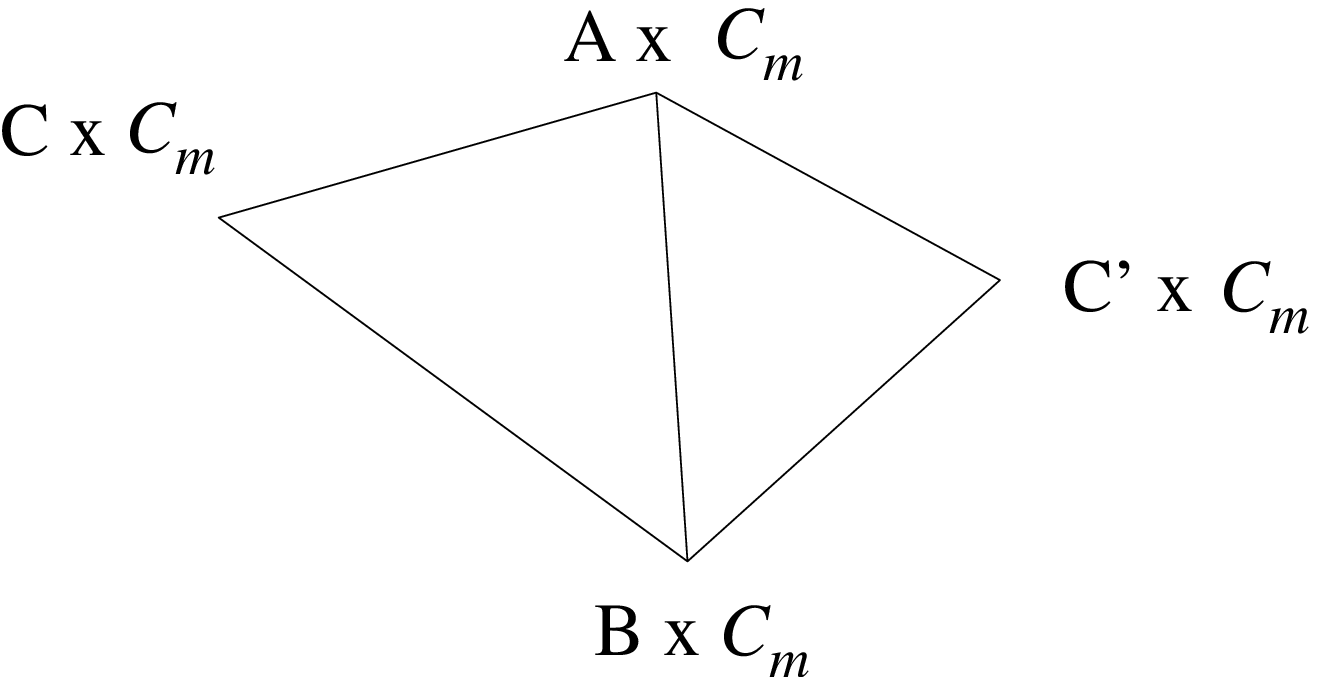}
        \caption{ Gluing triangulations in several copies of $C_3\times C_m$.}
                \label{fig:Cn}
                \label{fig:quadrangle}
        \end{figure}

\begin{thm}
If $m$ and $n$ are both even, then $P(m,n)$ can be triangulated with
$\frac{7}{2}mn -6(m+n) +8$ simplices.
\end{thm}

Observe that this coincides with the empirical values Table \ref{optitriang},
except for (6,6), (6,8) and (8,6), where it is two units above the value in the table.

\begin{proof}
We start with the triangulation $K$ obtained by replicating 
$n-2$ times the triangulation of $P(3,m)$ with 
$\frac{9}{2}-8$ simplices constructed in the previous theorem.
It will be important later that the triangulation (and labeling)
we choose for $C_n$ is exactly the one shown in Figure~\ref{fig:Cn},
with $n/2-2$ interior edges labeled $AB$, $n/2-1$ interior edges
labeled $BC$, and no interior edge labeled $AC$.

The proof consists on repeatedly using
the following trick: let us concentrate on two triangles of $C_n$
glued along a prism labeled, say, $AB$. We denote $C$ and $C'$
the vertices of $C_n$ opposite to the particular edge $AB$ 
we are considering (see again Figure~\ref{fig:quadrangle}).

Suppose there is a convex sub-polytope of the common $AB$-prism
that is triangulated in $K$. Suppose also
that all its simplices are joined to the same vertex $(C,i)$ in
        $C\times C_m$ (hence also to $(C',i)$ in $C'\times C_m$). Then, we have a
        bipyramid (a suspension) over $Q$, let us call it $SQ$,
        triangulated by first
        decomposing it into its two pyramids. One would expect
        that a more efficient way of triangulating $SQ$ is to join its axis
        to all the ``equatorial'' boundary simplices (that is to say,
        to all the triangles in $\partial Q$).

Things are actually a bit more complicated than suggested by the above
sentence. In the sentence we are implicitly assuming that the segment 
joining the two apices of the pyramids intersects the interior of $Q$
(otherwise we do not have a geometric bipyramid).  But it is not easy to guarantee that
this is indeed the case and, moreover, it is not
the most efficient way of doing things.

Indeed, if the axis intersects the interior of $Q$, then the number of
simplices that we get when we retriangulate equals the number of
triangles in $\partial Q$. But suppose, instead, that the axis
intersects a boundary point $x$ of $Q$. Then, we can retriangulate
by joining the axis to the triangles in facets of $\partial Q$ that
do not contain $x$. One problem with this is that then we have to take
care that the retriangulation of $SQ$ matches the rest of the
triangulation of $C_n\times C_m$ that we had. The way we guarantee
this is as follows: $Q$ is going to contain the segment
$[(A,i),(B,i)]$, and all the boundary faces of $Q$ containing $x$ are
going to be boundary faces of $AB\times C_m$ as well.
In the Cayley picture of Figure~\ref{fig:12gon-ABC} this property
corresponds to $Q$ containing a vertex of Minkowski sum and
part of the two boundary segments incident to it.
The consequence of this is that the segment $[(C,i),(C',i)]$
intersects $Q$ in a relative interior point of the edge
$[(A,i),(B,i)]$.
We will call that edge the 
\emph{distinguished edge} in the following discussion.

Figures~\ref{fig:12gon-AB} and~\ref{fig:12gon-BC} show how we implement this idea
in the $AB$ and $BC$ prisms, respectively. The shaded 
areas are the polygons $Q$ that we take in each of the prisms. 

       \begin{figure}[htb]
                \includegraphics[height=3in]{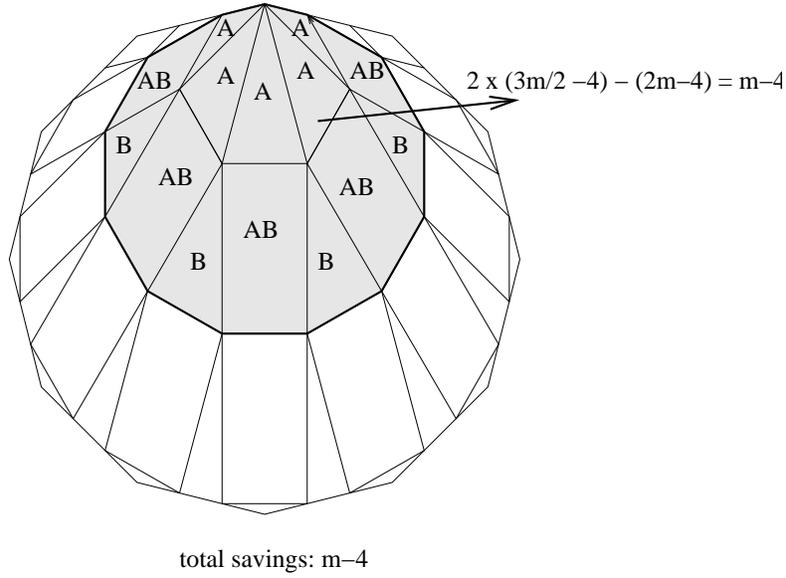}
        \caption{\small The region $Q$ in the facet $AB$.}
                \label{fig:12gon-AB}
        \end{figure}

In the $AB$-prism, the region $Q$ we consider is the
 the antiprism
        obtained from $AB\times C_m$ by cutting alternate corners
(this appears as a regular 12-gon in Figure~\ref{fig:12gon-AB})
together with one corner tetrahedron of the prism (the
small triangle in the top of the figure).
With that small triangle included, $Q$ is triangulated into $3m/2-4$
        simplices, so $SQ$ is triangulated into $3m-8$ simplices. The
        boundary of $Q$ has $2m-2$ triangles (two more than the
        antiprism would have), but two of them are incident to the
        distinguished edge. Hence, we can retriangulate $SQ$ into
        $2m-2$ simplices, saving us $m-4$ simplices in total.
Since we have  $n/2-2$ edges of type $AB$ we save $(n/2-2)(m-4)$ simplices.

In each prism $BC \times C_n$, we take several different polytopes $Q$ to apply the trick,
as shown in Figure~\ref{fig:12gon-BC}.
         \begin{figure}[htb]
                \includegraphics[height=3in]{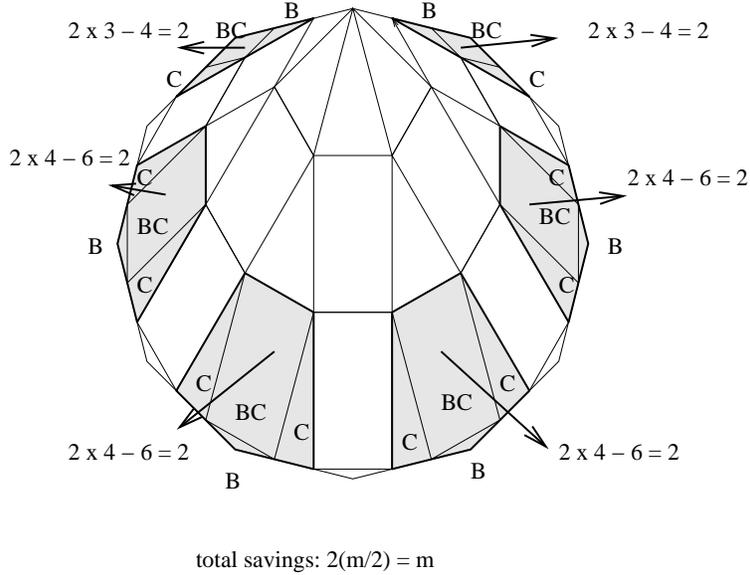}
        \caption{\small The regions $Q$ in the facet $BC$.}
                \label{fig:12gon-BC}
        \end{figure}
There are $m/2-2$ of a certain type and two of another. The type of each
is very easy to deduce from the figure: the two special ones are
triangular prisms (Cayley embedding of two equal triangles, one of
type $B$ and one of type $C$) and the other ones are 3-dimensional cubes with one
vertex truncated (Cayley embedding of a quadrilateral and a triangle
made with three of its four vertices. In the first type the original
triangulation of $SQ$ has 6 simplices, and we substitute them by four
simplices: the distinguished edge of the triangular prism is one
connecting the two opposite triangles, so there are four triangles in
facets not containing the point $x$. In the other type we originally have a
triangulation into 8 simplices and substitute it by one with only 6:
the distinguished edge is incident to two quadrilateral facets of $Q$,
and there remain another quadrilateral one (two triangles) 
plus four triangular ones. In total, we are decreasing the number of
simplices by $m$. Since we have $n/2-1$ prisms of this type,
we save $m(n/2-1)$ simplices in total.

Summing up, our final triangulation has
\[
(n-2)(9m/2 -8) - (n/2-2)(m-4) - (n/2-1)m =
7nm/2 -6n -6m + 8
\]
simplices, as claimed.

There is still one more thing that needs to be said in order to
justify correctness of the construction. In the triangulation of 
(almost all of) each copy of $C_3\times C_n$ we have done
changes to some pyramids with base on the $AB$ side and some on the
$BC$ side. Of course, for this to be possible we need these pyramids
to be disjoint. That they indeed are disjoint is easy to check in Figures
\ref{fig:12gon-AB} and \ref{fig:12gon-BC}. It just amounts to observing
that the shaded regions in the two pictures do not overlap.
\end{proof}

To finish the proof of Theorem~\ref{triangulations}, only the equations for $P(4,m)$ 
remain. The upper bound is just the substitution of $n=4$ in 
equation~(\ref{sevenhalves}). The idea for the lower bound is similar to the one in
Theorem~\ref{prism}.

\begin{thm}
\label{4xm}
The polytopal Gromov norm of $P(4,m)$ is at least $3\lceil 5(m-2)/2\rceil$.
\end{thm}

\begin{proof}
We regard $P(4,m)$ as a prism over the prism over an $m$-gon. That is to say, we regard its vertices
as lying in a ``bottom prism'' and a ``top prism''. Let $\alpha$ be an affine simplicial chain representing the
top relative homology class. As usual, we assume that the vertices of $\alpha$ are vertices of $P(4,m)$ and
that no four simplex in $\alpha$ lies in the boundary. Then, the simplices in $\alpha$ are of four types, depending on the number of
vertices they have on the top prism: we call them ``bottom'', ``half-bottom'', ``half-top'' and ``top'' simplices. The bottom and
top simplices need to cover the bottom and top prisms. By Theorem \ref{prism} there are at least $\lceil 5(m-2)/2\rceil$
of each type in $\alpha$. Also, by the same argument as in the proof of Theorem \ref{prism}, the numbers of half-bottom
and half-top simplices are each equal to at least that number, giving the total of $3\lceil 5(m-2)/2\rceil$ (rounded up to an even number).
\end{proof}

\section{A Binary Cover of $P(m,m)$ with $\frac{13m^2}{4}-\frac{19m}{2}$ Simplices} 
\label{upperbound}

Clearly, the definition of $||P||$ allows for much more freedom than 
using triangulations of $P$, in order to get upper bounds. Here we use
binary covers of $P$ for this purpose.

Recall that a pseudo-manifold is a simplicial complex of pure dimension
in which every codimension-one simplex lies in at most two
full-dimensional ones. Its boundary consists of the codimension-one
simplices that lie only in one full-dimensional simplex.
A \emph{binary
cover} of an $n$-dimensional polytope $P$ is a continuous 
 map $f:K\to P$ from  an oriented pseudo-manifold $K$ of dimension $n$
with the property that $f$ is linear on every
simplex and it restricts to a degree 1 map from
$\partial K$ to $\partial P$.

\begin{rem}
Every binary cover can be
homotoped to one that sends vertices of $K$ to vertices of $P$. Just choose,
for each vertex $v$ of $K$, a vertex of the minimal face of $P$ containing $f(v)$.
\end{rem}

\begin{lem}
\label{lem:binary}
If $f:K \to P$ is a binary cover of the polytope $P$, then $||P||$ 
is at most equal to the number of full-dimensional simplices in $K$.
\end{lem}

\begin{proof}
Since $K$ is a simplicial complex, there is an obvious chain
associated to it in which every top-dimensional simplex has weight
$1$ (the fact that $K$ is oriented is important here). 
We denote this chain by $K$ as well. The induced chain $f_*(K)$
is an affine chain of the polytope $P$. Because every codimension-one
simplex of $K$ lies in at most two full-dimensional ones, $f_*(K)$ is
a cycle in $S(P,\partial P)$. Because $f$  restricted to the boundary has
degree $1$ it follows (via Mayer-Vietoris) that $f$ itself
has degree $1$,
so $f_*(K)$ represents the fundamental class
$[P,\partial P]$. Therefore, $||P||$ is at most equal to the number of
simplices of $K$. 
\end{proof}

\bigskip

In this section, we exhibit two binary covers of $P(m,m)$,
for $m$ even. One has $\frac{13m^2}{4}$ and the other one
slightly less. Instead of describing the pseudo-manifold $K$,
we list the images of its simplices in $P(m,m)$. The pseudo-manifold
structure will be discussed later. We label the vertices of $P(m,n)$ by $(i,j)$ for $i,j=1,\dots,m$,
in the obvious way. Indices are regarded modulo $m$, and to list each simplex we give
its vertices.
The first list is:

\begin{figure}
  \centerline{
   \includegraphics[height=7 in]{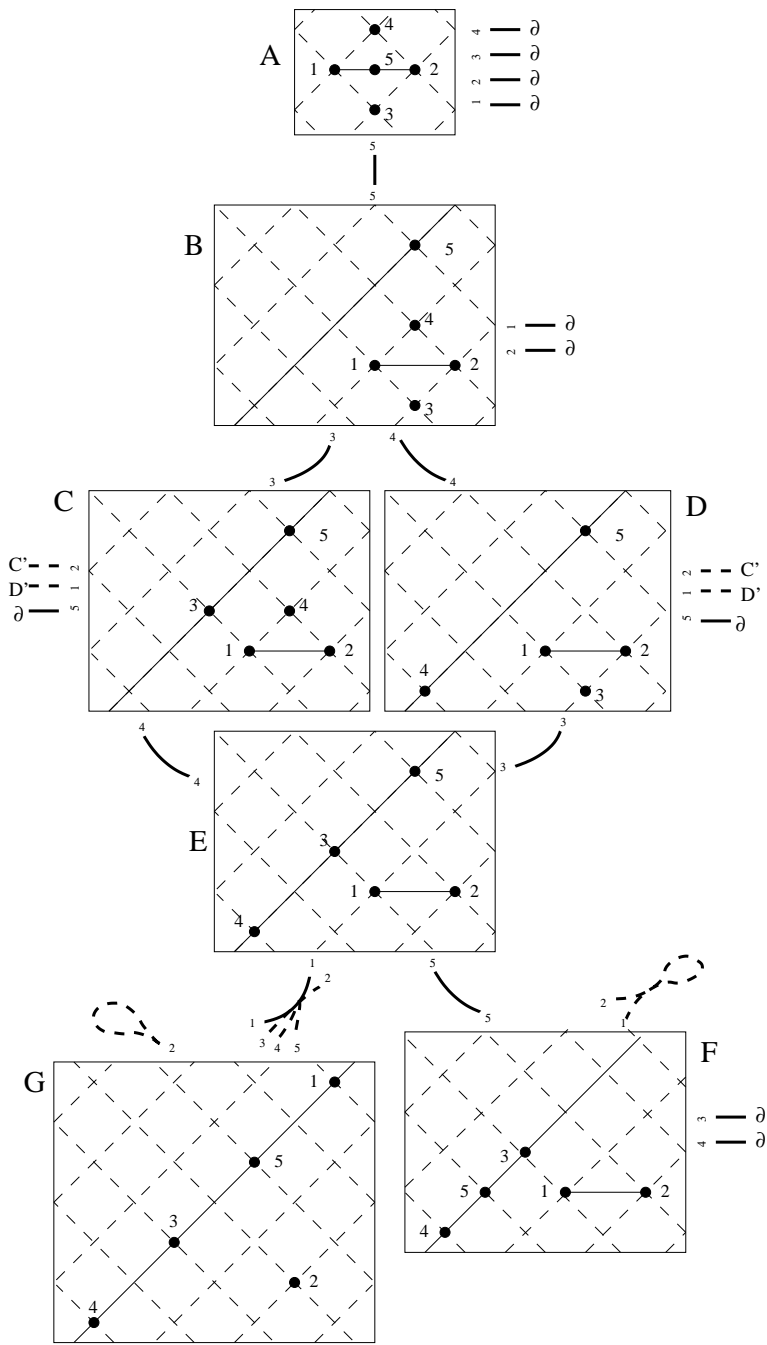}
                }
   \caption{ A small binary cover of $P(m,m)$.}
   \label{binary}
\end{figure}

\begin{enumerate}
\item[(1)]
 For each of the $m^2/4$ values of $(i,j)$ with $i$ even and $j$ odd,
the following six simplices, all of which contain the vertices
$(i-1,j)$ and $(i+1,j)$:
\smallskip

\begin{enumerate}
\item[(A)] The corner simplex at $(i,j)$. 
A corner simplex consists of $(i,j)$ and its four neighbors
$(i-1,j)$, $(i+1,j)$, $(i,j-1)$, $(i,j+1)$.

\item[(B)] The simplex 
$(i-1,j)$, $(i+1,j)$, $(i,j-1)$, $(i,j+1)$, $(i,i)$. 

\item[(C)] The simplex 
$(i-1,j)$, $(i+1,j)$, $(j+1,j+1)$, $(i,j+1)$, $(i,i)$. 

\item[(D)] The simplex 
$(i-1,j)$, $(i+1,j)$, $(i,j-1)$, $(j-1,j-1)$, $(i,i)$. 

\item[(E)] The simplex 
$(i-1,j)$, $(i+1,j)$, $(j+1,j+1)$, $(j-1,j-1)$, $(i,i)$. 

\item[(F)] The simplex 
$(i-1,j)$, $(i+1,j)$, $(j+1,j+1)$, $(j-1,j-1)$, $(j,j)$. 
\end{enumerate}

\item[(2)] Symmetrically, for each of the $m^2/4$ values of $(j,i)$ with $j$ odd and $i$ even,
the six simplices (A'), (B'), (C'), (D'), (E') and (F') obtained from the previous six
by exchanging $i$ and $j$.

\item[(3)] Finally, for each of the $m^2/4$ values of $(i,j)$ with $i$ and $j$ odd,
the simplex $(i,j)$, $(i-1,i-1)$, $(i+1,i+1)$, $(j-1,j-1)$, $(j+1,j+1)$.

\end{enumerate}

This gives 13 types of simplices, which we will refer to as
(A), (B), (C), (D), (E), (F), (A'), (B'), (C'), (D'), (E'), (F') and (G).
There are $m^2/4$ of each type.
Figure~\ref{binary} schematically shows the construction. It depicts
one simplex of each of the types (A) through (G), 
 each drawn as a set of five points in a 
2-dimensional grid. The dashed diagonal line in six of the simplices
represents the set of vertices $(i,i)$ of $P(m,m)$. 

Observe that some of the simplices are degenerate (they are not full-dimensional
or they even have repeated vertices). This happens, for example, for the simplices
(B) through (F) if $i-j=\pm 1$, or for any simplex $G$ if the factor polygons
of $P(m,m)$ are equal (the diagonal of the grid represents then a 2-plane
in $\R^4$, and the ``4-simplex'' $G$ has a 3-face lying in it).

To help check that this is indeed a binary cover, the incidences between simplices
are marked in the figure. More precisely, each simplex has five ``bonds'' to either other simplices,
or to the symbol $\partial$, representing the boundary of $K$.
The five vertices in each simplex are labeled 1 through 5 and the bond labeled $i$ on one side
represents the facet opposite to vertex $i$ on that simplex. A bond is drawn solid if it
joins exactly the simplices in the figure (or if the corresponding facet lies in the boundary of $P(m,m)$
and is drawn dashed if it is between the simplex in the picture and one not in the picture, (but of the
same type). It is left to the reader to check that,
with the glueings specified by the bonds in the picture, the list of simplices
is indeed an oriented pseudo-manifold with boundary.

Now we try to understand how the boundary of the pseudo-manifold covers the boundary of $P(m,m)$. 
The first check is that, indeed, all the facets of simplices with bonds to the symbol $\partial$
lie in the boundary of $P(m,m)$. Next, we concentrate on a facet of $P(m,m)$, say the prism consisting on the 
vertices $(i,*)$ and $(i+1,*)$ for some (say, even) $i$. Each simplex of type (A) or (A') ``centered'' at a point
$(i,j)$ or $(i+1,j)$ contains a facet on our prism, and it cuts a corner of it. After all of them are removed, what remains
is an $m/2$-antiprism consisting of the vertices $(i,j)$ for even $j$ and $(i+1,j)$ for odd $j$.
The other simplices with facets in our prism are those of types (B'), (C) and (F). This is so because 
(B), (C'), (D') and (F') only contain facets on ``vertical'' prisms, and (D) contains a facet in every other horizontal prism, but not
the one we are considering. It can be easily checked that the facets that (B'), (C) and (F) have in our antiprism produce the following
degree one cover of it: consider the cover of the $m/2$-gon $(i,*)$ (where ``*'' is meant to be even) obtained by coning 
$(i,i)$ to the boundary. Then join this cover, as well as the $m$ triangular faces of the antiprism, to $(i+1,i+1)$.

With all this we conclude that:

\begin{thm}
The above list of $13m^2/4$ simplices forms a binary cover of $P(m,m)$.
\end{thm}

Our next goal is to show that this binary cover, call it $\alpha$, contains as a proper 
subset an even smaller binary cover. This is obtained by deleting from the initial 
binary cover all the simplices with repeated vertices (but this condition is not enough to guarantee
that they can be removed. The reader should check that after the removal, and with some
minor regluing, we still have an oriented  pseudo-manifold):

\begin{enumerate}

\item The simplices of types (B) and (B') for which $i-j=\pm 1$ ($2m$ of them).
\item The simplices of types (C) and C') for which $i-j=1$ ($m$ of them).
\item The simplices of types (D) and (D') for which $i-j=- 1$ ($m$ of them).
\item The simplices of types (E) and (E') for which $i-j=\pm 1$ ($2m$ of them).
\item The simplices of types (F) and (F') for which $i-j=\pm 1$ ($2m$ of them).
\item The simplices of type (G) and (E') for which $i-j=0$ or $\pm2$ ($3m/2$ of them, if $m\ge 6$).

\end{enumerate}

This deletes $19m/2$ simplices from the initial list. We leave it to the reader to
check that indeed this is a binary cover.

\begin{cor}
For every even $m\ge 6$, $P(m,m)$ has a binary cover with
$13m^2/4-19m/2$ simplices.
\end{cor}

\begin{rem}
Observe that this new binary cover still has some degenerate simplices,
at least if we assume the two $C_m$ factors in $P(m,m)$ to be equal. 
For example, the $m^2/4-3m/2$ simplices of type $G$ all have a 3-face
contained in a 2-plane.
Even though
they do not cover any ``space'', their removal would leave
some interior tetrahedra unmatched. In other words, the $3m^2-8m$
simplices of  types $A$ through $F'$ form a
cover of $P(m,m)$ without overlaps, but this cover is
insufficient to make a statement about the Gromov norm because
some faces are unmatched.
\end{rem}

\section{A Lower Bound for the polytopal Gromov norm} \label{lowerbound}

In this section we prove a lower bound for the polytopal Gromov
norm of $P$. by counting (with weights) certain incidences in affine chains
of $S(P, \partial P)$. 

Each affine 4-simplex $\sigma \in P(m,n)$ has 20
triangle-tetrahedron incidences. We say that one of these incidences is a {\it titap}
incidence if the tetrahedron is contained in a facet (prism) of $P(m,n)$ and the triangle
is interior to that facet. (``Titap'' is short for ``\emph{triangle interior to a prism}''). We denote the
number of titap incidences in $\sigma$ as $\titap(\sigma)$. Similarly, for an affine chain
 \begin{displaymath}
c= \Sigma_i \, w_i \sigma_i \in S(P,\partial P)
\end{displaymath}
we define
 \begin{displaymath}
\titap(c)= \Sigma_i \, |w_i| \titap( \sigma_i).
\end{displaymath}

\begin{lem}
\label{titap-prism}
For every affine chain $c \in S(P,\partial P)$,
\[
\titap(c) \ge 12mn - 16m -16n.
\]
\end{lem}

\begin{proof}
Clearly, the titap incidences in $c$ can be counted by adding the titap incidences in the restrictions of
$c$ to the individual boundary prisms of $P(m,n)$ (because each titap incidence belongs to one and only one prism).

Let $c$ be an affine chain and let $c'$ be its restriction to a certain $m$-prism. As in 
the proof of Theorem \ref{prism}, we classify the tetrahedra in $c'$ as ``bottom'', ``middle''
and ``top'', depending on their number of vertices in the bottom and top $m$-gons of the prism.
We count
titap incidences in the three groups of tetrahedra separately.

Each bottom tetrahedron $\tau$ has a unique triangle $\rho$ in the bottom $m$-gon of the prism. Clearly, if
an edge of $\rho$ is interior to the $m$-gon, the corresponding triangle in $\tau$
is a titap incidence. Since the bottom triangles must cover the bottom $m$-gon, they produce at
least $2(m-3)$ of these incidences (because a binary cover of the bottom $m$-gon has at least $m-3$
interior edges, each in at least two triangles). Similarly, there are at least $2(m-3)$ titap incidences in
top tetrahedra.

Each middle tetrahedron has two bottom triangles (with 2 vertices in the bottom $m$-gon and one
in the top $m$-gon) and two top triangles. Some of these triangles may be in vertical faces of the prism, but
(as in the proof of Theorem \ref{prism}) we at least know that the bottom triangles cover
the $m$-gon, when projected to it, and the same for the top triangles. Hence, there are at least 
$m-2$ of each type that are not vertical. Hence, middle tetrahedra 
produce at least $2m-4$ titap incidences.
Adding this up, we conclude that an $m$-prism contains at least
$
2(m-3) + 2(m-3) + 2m-4  = 6m-16
$
titap incidences. Adding over the $n$ $m$-prisms plus $m$ $n$-prisms gives the statement.
\end{proof}

\begin{lem}
\label{titap-simplex}
Let $\sigma$ be an affine simplex in $P(m,n)$, not contained in the boundary.
Then,
$
\titap(\sigma)\le 6.
$
\end{lem}

\begin{proof}
Let $k$ be the number of facets of $\sigma$ that lie in the boundary of $P(m,n)$.
Since $\sigma$ is not contained in the boundary, the $k$ boundary tetrahedra in $\sigma$
lie each in a different facet of $P(m,n)$. In particular, the common triangle to two
of them is not interior to a prism, and does not produce a titap incidence.
Then, each of the $k$ tetrahedra produces at most $4-(k-1)=5-k$ titap incidences,
because $k-1$ of its four triangles are used in adjacencies to other boundary tetrahedra.
Hence,  $\sigma$ has at most $k(5-k)$ titap incidences.
The maximum of $k(k-1)$ is 6, achieved for $k=2$ or 3.
\end{proof}

\begin{cor} 
$
T(m,n) \ge ||P(m,n)|| \ge 2mn-8(m+n)/3.
$
\end{cor}

\begin{proof}
For every affine chain $c= \Sigma_i \, w_i \sigma_i \in S(P,\partial P)$,
\[
12mn-16m-16n \le \titap(c) =\Sigma_i \, |w_i| \titap(\sigma_i) \le 6 \Sigma_i \, |w_i| = 6||c||,
\]
where the two inequalities come from the previous two lemmas.
\end{proof}

\begin{rem}
We do not believe our lower bound to be very close to the real value of $||P(m,n)||$ or
$T(m,n)$, because it is based in a very specific type of incidence. Our conjecture is that 
$||P(m,n)||$ is closer to the upper bound obtained in Section \ref{upperbound}, perhaps in
$3mn \pm O(m+n)$.

Observe also that our lower bound can be slightly improved if we restrict our attention to
corner-cutting triangulations, that is to say, triangulations that first cut $mn/2$
vertices of $P(m,n)$ via corner 4-simplices, for $m$ and $n$  even.
The $mn/2$ corner simplices produce only $2mn$ titap incidences, and we need at least
another $(10mn-16n-16m)/6$ simplices to produce the rest, giving a total of
at least $13mn/6 -O(m+n)$ simplices.
\end{rem}

%
%
%
%
%

\noindent{\bf Acknowledgements}
We are grateful to Greg Kuperberg who first proposed us to study this problem.
We also received several useful comments from Bernd
Sturmfels and G\"unter Ziegler. The last three authors gratefully
acknowledge support of the Mathematical Science Research Institute at
UC Berkeley where most of this research was conducted. The first
and second author were partially supported by NSF VIGRE Grant
DMS-0135345 and the second author gratefully acknowledges support
from NSF grant DMS-0309694. Work of the fourth author is 
supported by grant BMF2001-1153 of Spanish Direcci\'on
General de Investigaci\'on Cient\'{\i}fica.

\clearpage

\section*{ERRATUM}

Michelle Bucher-Karlsson has communicated to us that our proof of
Lemma 2.7 is wrong. We say:

\begin{quote}
\emph{
For every $h>0$, there is a regular $4h$-gon $F_h$ with
all interior angles equal to $2\pi/4h$. We choose $h$ large enough so
that every pair of nonadjacent sides of $F_h$ is at least a distance
$d$ apart.
}
\end{quote}

In this sentence, $d$ could be arbitrarily large, since it was defined as the maximum diameter of an (a-priori, arbitrary) simplex in $\H^2 \times \H^2$.
Our choice of $h$ is then in contradiction with the following statement,
communicated to us by Bucher-Karlsson:

\bigskip
\noindent
\textbf{Proposition.}
\emph{
For any $h>0$, the distance between the midpoints of two
adjacent edges in $F_h$ is smaller than $\operatorname{arccosh}(3)\sim 1.763$. 
In particular, $F_h$ contains two non-adjacent edges 
at distance smaller than $2 \operatorname{arccosh}(3)$.
}

\begin{proof}
Given a hyperbolic geodesic triangle with angles $\alpha$, $\beta$,
$\gamma$ and opposite sides of lengths $a$, $b$, $c$ respectively, the second
cosine rule for hyperbolic triangles states that
\[
    \sin(\beta)\sin(\gamma)\cosh(a)=\cos(\alpha)+\cos(\beta)\cos(\gamma).
\]

Consider the geodesic triangle with vertices the midpoints of two
adjacent edges and the center of $F_h$. The angle at the center is equal
to $2\pi/4h$ and, by a symmetry argument,  
the angles at the two other
corners are both equal to $\pi/4$. Thus, by the second cosine rule 
 the distance between the midpoints of two adjacent edges
is equal to
\[
    \operatorname{arccosh}\left(2\cos\left(2\pi/4h\right)+1\right),
\]
and is hence bounded by $\operatorname{arccosh}(3)$.
\end{proof}

Without Lemma~2.7, one direction of the equality in our Theorem~2.3,
and in its generalization Theorem~1.1, is invalid. The correct
statements must now be:
\bigskip

\setcounter{section}{1}
\setcounter{pro}{0}

\begin{thm}
\label{surfaces}
Let $||P||$ denote the polytopal Gromov norm of a polytope $P$.
Then, the Gromov norm of the product $\Sigma_g\times\Sigma_h$
of two surfaces of genera $g$ and $h$ equals
\begin{equation}
\frac{||\Sigma_g\times\Sigma_h||}{(g-1)(h-1)} \le 16\lim_{n,m\to\infty}\frac{||P(n,m)||}{nm}= 16\inf_{n,m}\frac{||P(n,m)||}{nm}. 
\label{eqn-covering}
\end{equation}
\end{thm}

\setcounter{section}{2}
\setcounter{pro}{2}

\begin{thm}
\begin{equation}
||\Sigma_2 \times \Sigma_2|| \le \lim_{n,m \to \infty} \,  \frac{16||P(n,m)||}{(n-2)(m-2)}.
\end{equation}
\end{thm}

The lower bound $32(g-1)(h-1) \le ||\Sigma_g\times \Sigma_h||$ that we
gave in Corollary~1.5 is also invalid. In fact,
Bucher-Karlsson has computed \emph{exactly} the value of
$\Sigma_g\times \Sigma_h$:

\medskip
\noindent
\textbf{Theorem.} (Bucher-Karlsson~\cite{BK1})
\emph{Let $M$ be a closed, oriented Riemannian manifold whose universal
cover is isometric to $\H_2 \times \H_2$ . Then
\[
||M||=6 \chi(M).
\]
}

As an update, in a more recent paper~\cite{BK2} the same author has
improved our lower bound for the polytopal Gromov norm of the product
of two polygons.  We proved $||P(m,n)|| \ge 2mn-O(m+n)$ and she gets
$||P(m,n)|| \ge 3.125 mn - 5(m+n) +6$. 

This confirms what we said in Remark 5.4: ``We do not believe our
lower bound to be very close to the real value of $||P(m,n)||$ or
$T(m,n)$. Our conjecture is that $||P(m,n)||$ is closer to the upper
bound\footnote{$||P(m,m)|| \le 3.25 m^2 - \Omega(m+n)$} obtained in Section
4''. 
We were only wrong in our final guess ``perhaps in $3mn \pm
O(m+n)$''.

\end{document}